\documentclass[a4paper,10pt]{article}

\usepackage{amsfonts}
\usepackage{latexsym}
\usepackage{epsfig}

\usepackage{amssymb}

\usepackage{graphicx}
\usepackage{oldgerm}
\usepackage{theorem}

\def\ZV{{\mathbb{Z}}^{N}_{<}}
\def\RV{{\mathbb{R}}^{N}_{<}}
\def\S{{\bf S}}
\def\u{{\bf u}}
\def\v{{\bf v}}
\def\x{{\bf x}}
\def\y{{\bf y}}
\def\z{{\bf z}}

\def\0{{\bf 0}}

\def\X{{\bf X}}
\def\Y{{\bf Y}}

\def\cN{{\cal N}_N}

\def\P{{\bf P}}

\def\vlambda{\mib{\lambda}}
\def\vxi{\mib{\xi}}

\def\Z{\mathbb{Z}}
\def\R{\mathbb{R}}
\def\C{\mathbb{C}}
\def\N{\mathbb{N}}

\def\Sym{\mathfrak{S}}

\usepackage{theorem}

\newcommand{\mib}[1]{\mbox{\boldmath $#1$}}

\theorembodyfont{\itshape}
\newtheorem{thm}{Theorem}
\newtheorem{lem}[thm]{Lemma}
\newtheorem{cor}[thm]{Corollary}
\newtheorem{prop}[thm]{Proposition}

\newcommand{\qed}{\hbox{\rule[-2pt]{3pt}{6pt}}}
\begin{document}
\pagestyle{plain}
\vskip 0.3cm

\begin{center}
{\bf \LARGE{Nonintersecting Paths, Noncolliding Diffusion Processes}}
\vskip 0.3cm

{\bf \LARGE{and Representation Theory}}
\footnote{This manuscript is based on the talks 
at the Institute of Applied Mathematics, Chinese Academy of 
Science, Beijing, China, on 27th September 2004,
and in the workshop 
{\it `Combinatorial Methods in Representation Theory and
Their Applications'} (19-22 October 2004) 
at the Research Institute for Mathematical Sciences,
Kyoto University, Kyoto, Japan; 
to be published in {\it RIMS Kokyuroku}.}

\vskip 0.5cm

\begin{large}
{\sf Makoto Katori}
\end{large}
\footnote{\sf Electronic mail: katori@phys.chuo-u.ac.jp}\\
\vskip 0.3cm
{\it Department of Physics,
Faculty of Science and Engineering,\\
Chuo University, 
Kasuga, Bunkyo-ku, 
Tokyo 112-8551, Japan}

\vskip 0.3cm

\begin{large}
{\sf Hideki Tanemura}
\end{large}
\footnote{\sf Electronic mail: tanemura@math.s.chiba-u.ac.jp}\\
\vskip 0.3cm
{\it Department of Mathematics and Informatics,
Faculty of Science,\\
Chiba University, 
1-33 Yayoi-cho, Inage-ku,
Chiba 263-8522, Japan}
\end{center}

\vskip 0.5cm

\noindent {\bf Abstract} \quad
The system of one-dimensional symmetric simple random walks,
in which none of walkers have met others in a given time period,
is called the {\bf vicious walker model}. It was introduced by Michael
Fisher and applications of the model to various wetting and melting
phenomena were described in his Boltzmann medal lecture.
In the present report, we explain interesting connections
among representation theory, probability theory, 
and random matrix theory 
using this simple diffusion particle system.
Each vicious walk of $N$ walkers is represented by an $N$-tuple of
{\bf nonintersecting lattice paths} on the spatio-temporal plane.
There is established a simple bijection between
nonintersecting lattice paths and {\bf semistandard Young tableaux}.
Based on this bijection and some knowledge of symmetric
polynomials called the {\bf Schur functions}, we can give a
determinantal expression to the partition function of vicious
walks, which is regarded as a special case of the
{\bf Karlin-McGregor formula} in the probability theory
(or the {\bf Lindstr\"om-Gessel-Viennot formula} in the 
enumerative combinatorics). Due to a basic property of Schur
function, we can take the diffusion scaling limit of the vicious
walks and define a {\bf noncolliding system of Brownian particles}.
This diffusion process solves the stochastic differential
equations with the drift terms acting as the repulsive two-body
forces proportional to the inverse of distances between particles,
and thus it is identified with {\bf Dyson's Brownian motion model}.
In other words, the obtained noncolliding system of Brownian
particles is equivalent in distribution with the {\bf eigenvalue process
of a Hermitian matrix-valued process}.

\setcounter{equation}{0}
\section{Vicious Walks, Young Tableaux and
Schur Functions}

\begin{figure}[htpb]
\includegraphics[width=0.7\linewidth]{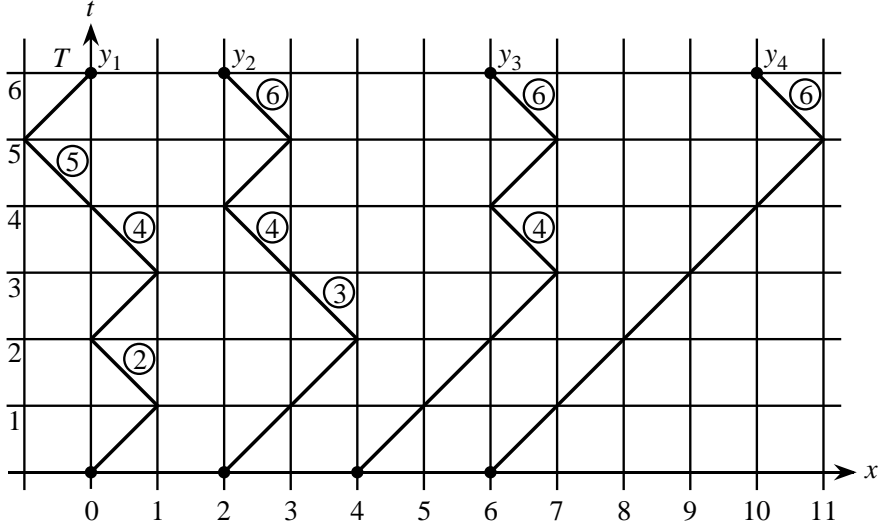}
\caption{An example of vicious walk 
in the case $N=4, T=6$. \label{fig:Fig1}}
\end{figure}

Let $\Big(\{\S(t)\}_{t=0,1,2,...},
{\sf P}^{\x} \Big)$ 
be the $N$-dimensional Markov chain
starting from\\
$\x=(x_1,x_2,\dots, x_N)$, 
such that the coordinates $S_{i}(t), i=1,2,\dots,N$,
are independent simple random walks on $\Z$.
We always take the starting point $\x$ from the set
$$
\ZV =
\Big\{ \x =(x_1,x_2,\dots,x_N) \in (2 \Z)^N \; :
x_{1}<x_{2}<...<x_{N} \Big\}.
$$
We consider the condition that any of walkers 
does not meet other walkers up to time $T>0$, 
{\it i.e.}
\begin{equation}
S_1(t) < S_2(t) < \cdots < S_N(t),
\quad t=1,2,...,T.
\label{eqn:nonint}
\end{equation}
We denote by ${\sf Q}^{\x}_T$ the conditional probability of 
${\sf P}^{\x}$
under the event 
$\Lambda_T = \Big\{S_1(t) < S_2(t) < \cdots < S_N(t), 
\; t=0,1,...,T \Big\}$.
M. Fisher called
the process $\Big(\{\S(t)\}_{t=0,1,2,...,T},
{\sf Q}^{\x}_T \Big)$ 
the {\bf vicious walker model} in his
Boltzmann medal lecture \cite{Fis84}.

We will assume the initial positions as
\begin{equation}
   S_{i}(0)=2(i-1), \qquad i=1,2, \cdots, N
\label{eqn:initial}
\end{equation}
in Sections 1 in this report.

Each realization of vicious walk is represented
by an $N$-tuple of {\bf nonintersecting
lattice paths} on the 1+1 spatio-temporal plane,
$\Z \times \{0,1, \cdots, T\}$.
An example is given by Figure 1 in the case that
four walkers ($N=4$) perform a noncolliding walk
up to time $T=6$.

\begin{figure}[htbp]
\begin{minipage}{0.5\linewidth}
\includegraphics[width=0.9\linewidth]{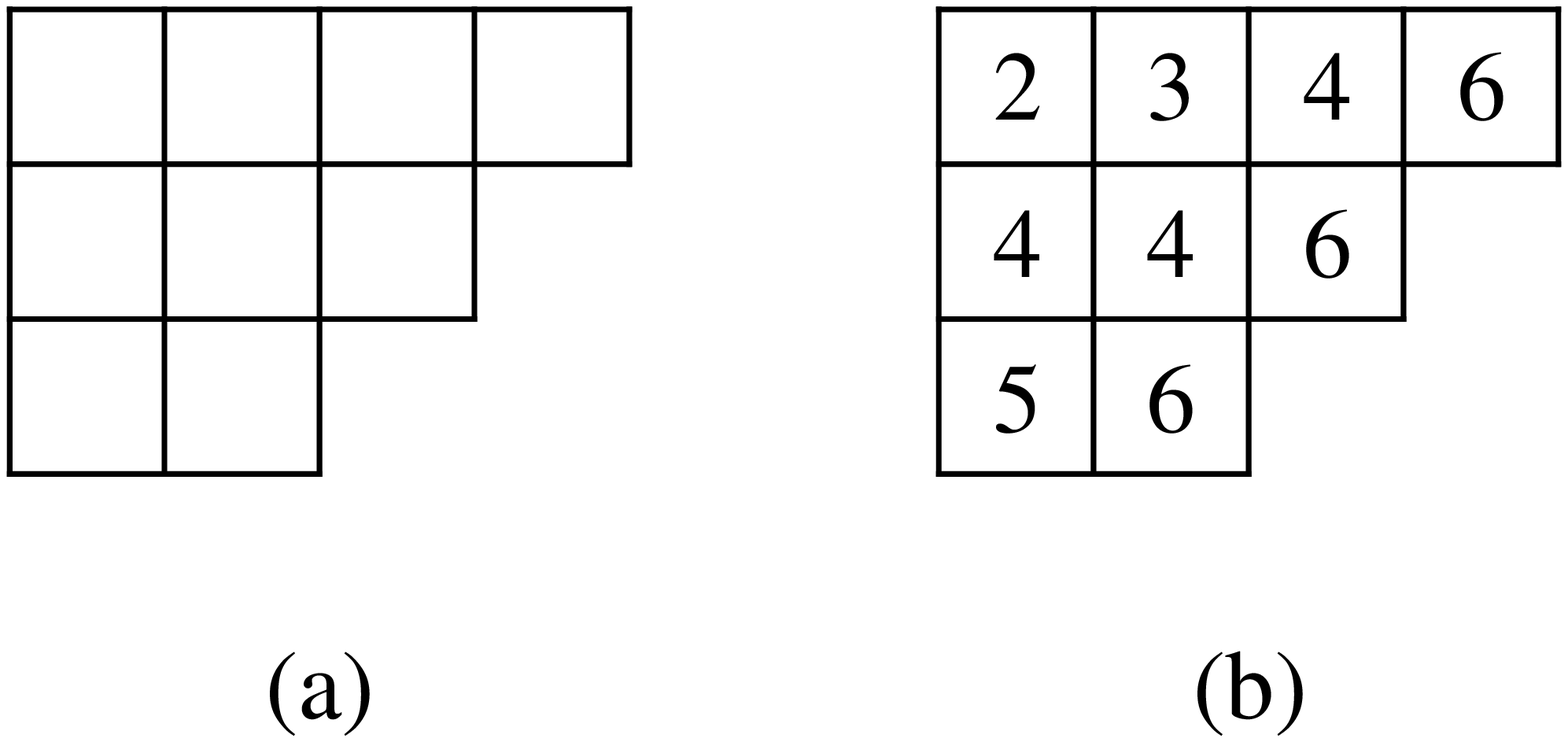}
\label{fig:Fig2}
\caption{(a) Young diagram and
(b) Young tableau ${\cal T}$ corresponding to
the vicious walk in Figure 1.}
\end{minipage}
\begin{minipage}{0.5\linewidth}
\includegraphics[width=0.9\linewidth]{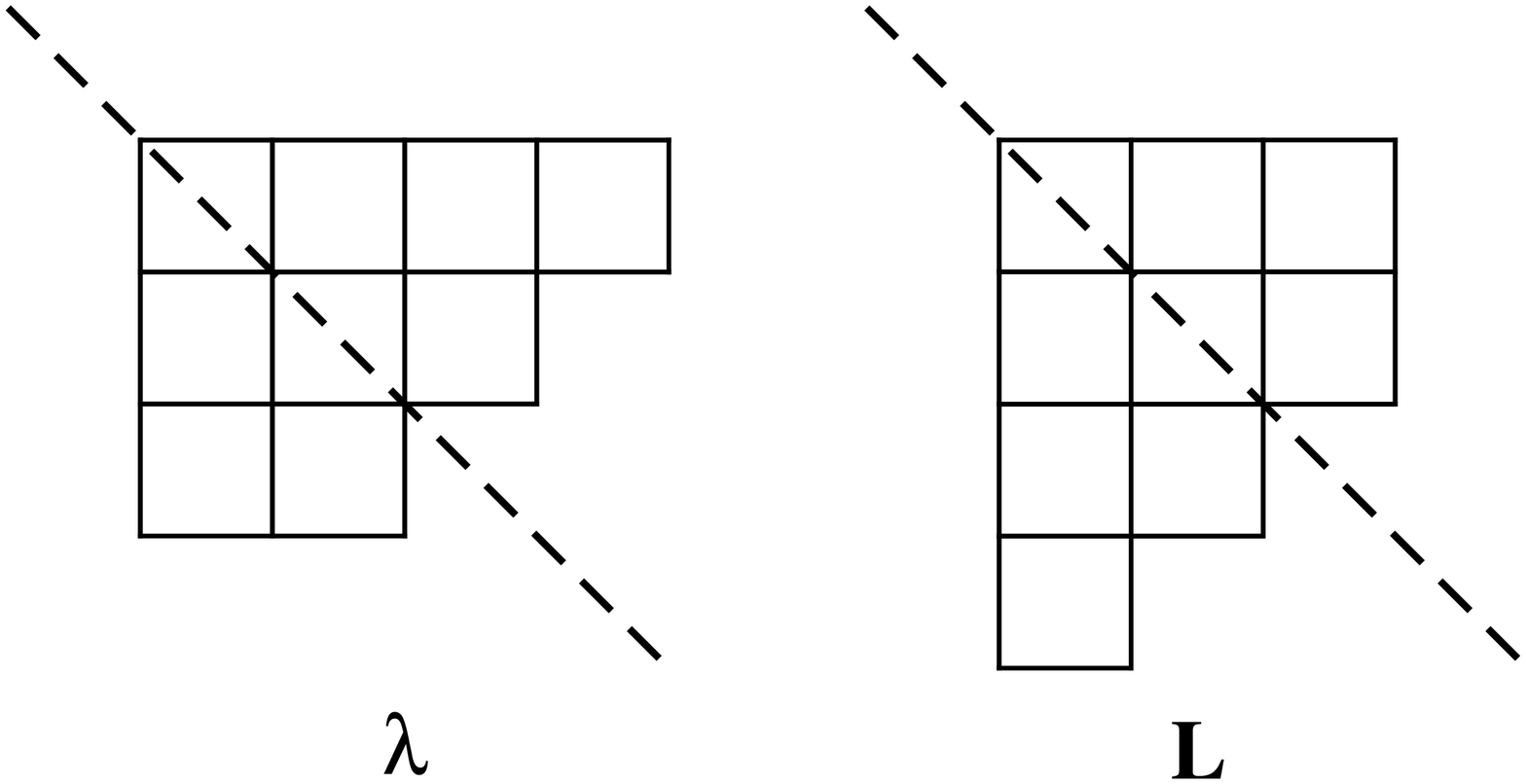}
\caption{A pair of conjugate YD's. }
\label{fig:Fig3}
\end{minipage}
\end{figure}

Bijection between such
nonintersecting lattice paths
and {\bf semistandard Young tableaux (SSYT)}, 
$
    {\cal T}=({\cal T}(i, j)),
$
is established by the following procedure
\cite{GOV98,KGV00}.
\begin{description}
\item{(1)} \quad For $1 \leq j \leq N$, let
$$
L_{j} \equiv \mbox{the number of leftward steps 
among $T$ steps of the $j$-th walker}.
$$
Draw a collection of boxes with $N$ columns,
in which the number of boxes in the $j$-th column
is $L_{j}$. 
(We number columns from the left to the right.)
Since ${\bf L} \equiv (L_{1}, L_{2}, L_{3}, L_{4})
=(3,3,2,1)$ in the walk shown in Figure 1,
we draw the collection of boxes as shown in Figure 2 (a) 
for this example.
\item{(2)} \quad For each walker, we label each leftward step
by the integer $\in \{1,2, \cdots, T\}$, which is the
time when that leftward step was done.
See Figure 1, in which labels of leftward steps
are indicated by integers in small circles 
associated with the line segments showing leftward steps.
Then for the $j$-th column of the collection of boxes,
fill the boxes by the labels of leftward steps
of the $j$-th walker, 
from the top to the bottom, $1 \leq j \leq N$.
For the walk given in Figure 1, we have the
boxes with integers shown in Figure 2 (b). Let
$$
{\cal T}(i, j)= \mbox{the integer in the box located in 
the $i$-th row and $j$-th column}. 
$$
For example, ${\cal T}(1,3)=4$ and ${\cal T}(3,1)=5$
in this case.
\end{description}
\begin{description}
\item[Remark 1.] \,
The above procedure with
the nonintersecting condition (\ref{eqn:nonint})
guarantees the inequalities
\begin{equation}
  L_{1} \geq L_{2} \geq \cdots \geq L_{N},
\label{eqn:Leq}
\end{equation}
and
\begin{eqnarray}
   && {\cal T}(i, j) < {\cal T}(i+1, j), \qquad
   \mbox{strictly increasing in each column},
   \nonumber\\
   && {\cal T}(i, j) \leq {\cal T}(i, j+1), \qquad
   \mbox{weakly increasing in each row}.
\label{eqn:Teq}
\end{eqnarray}
Assume that the number of rows in the collection of
boxes is $\ell$.
Let 
$$
\lambda_{i} \equiv \mbox{the number of boxes in the $i$-th
row}, \, i=1,2, \cdots, \ell.
$$
Then the inequalities (\ref{eqn:Leq}) imply
\begin{equation}
 \lambda_{1} \geq \lambda_{2} \geq \cdots 
 \geq \lambda_{\ell}.
\label{eqn:lameq}
\end{equation}
The collections of boxes with such conditions 
concerning the numbers of boxes in
rows (\ref{eqn:lameq}) 
(and in columns (\ref{eqn:Leq})) are called
{\bf Young diagram (YD)}.
The number of rows $\ell$ in the YD is
called the {\bf length} of YD. In the
present situation, $\ell \leq T$, in general.
(In our example in Figure 2 (b), 
$\ell=3$ for $T=6$.) 
YD's with integers with the conditions
(\ref{eqn:Teq}) are called 
{\bf semistandard Young tableaux (SSYT)}. 
\item[Remark 2.] \,
YD with $\lambda_{i}$ boxes in the $i$-th row,
$1 \leq i \leq \ell$,
is said to be the YD of {\bf shape} 
$\vlambda=(\lambda_{1}, \cdots, \lambda_{\ell})$.
The YD with the shape 
${\bf L}=(L_{1}, \cdots, L_{N})$ is regarded as
the {\bf conjugate} of the YD with the shape $\vlambda$
and denoted by
$$
  {\bf L}=\widetilde{\vlambda}.
\label{eqn:conjugate1}
$$
As shown in Figure 3, they are mirror images
with respect to the diagonal line.
\end{description}

A sequence of integers with the condition
(\ref{eqn:lameq}), that is,
$$
\vlambda=(\lambda_{1}, \lambda_{2}, \cdots, \lambda_{T})
\quad \mbox{with} \quad
\lambda_{i} \in \N \equiv \Big\{x \in \Z, x \geq 0 \Big\}, \,
1 \leq i \leq T,
\quad \lambda_{1} \geq \lambda_{2} \geq \cdots
\geq \lambda_{T},
$$
is regarded as a {\bf partition} of an integer
$n=\sum_{i=1}^{T} \lambda_{i}$.
We introduce a set of $T$ variables
$\z=\{z_{1}, z_{2}, \cdots, z_{T}\} \in \C^{T}$
and define a monomial
\begin{eqnarray}
  \z^{{\cal T}} &=& \prod_{(i, j)} z_{T(i, j)} \nonumber\\
  &=& \prod_{k=1}^{T} z_{k}^{\mbox{\footnotesize $\#$ of times that
  the integer $k$ occurs in ${\cal T}$}}.
\nonumber
\end{eqnarray}
For example, the monomial corresponding to the
SSYT ${\cal T}$ shown in Figure 2 (b) is
\begin{eqnarray}
\z^{{\cal T}} &=& z_{2} \times z_{3} \times z_{4} \times z_{6} 
\nonumber\\
&\times& z_{4} \times z_{4} \times z_{6} \nonumber\\
&\times& z_{5} \times z_{6} \nonumber\\
&=& z_{2} z_{3} z_{4}^{3} z_{5} z_{6}^3. \nonumber
\end{eqnarray}
Notice that for one YD with a given shape $\vlambda$,
there are different ways of filling
boxes with integers to make SSYT's satisfying 
the conditions (\ref{eqn:Teq}).
For each YD with shape $\vlambda$, we define
a polynomial of $\z=\{z_{1}, \cdots, z_{T}\}$ by
summing $\z^{{\cal T}}$ over all SSYT defined on the YD:
$$
s_{\lambda}(z_{1}, z_{2}, \cdots, z_{T})
= \sum_{{\cal T}: \mbox{{\scriptsize all SSYT with the same
shape $\vlambda$}}} \z^{{\cal T}}.
$$
This polynomial is called the
{\bf Schur function} indexed by
(the partition/YD with shape) $\vlambda$
on $\z=(z_{1}, \cdots, z_{T})$.
We can prove the two formulae
({\bf Jacobi-Trudi formulae}).
The first one is the following.
\vskip 0.5cm
\fbox{\parbox{16.5cm}{
\begin{lem}
\label{thm:JT1}
$$
s_{\lambda}(z_{1}, \cdots, z_{T})
= \frac{\displaystyle{ \det_{1 \leq i, j \leq T}
\Bigg[ z_{i}^{\lambda_{j}+T-j} \Bigg]}}
{\displaystyle{\det_{1 \leq i, j \leq T}
\Bigg[ z_{i}^{T-j} \Bigg] }},
$$
where the denominator is the {\bf Vandermonde determinant}
evaluated as the {\bf product of differences},
\begin{equation}
\det_{1 \leq i, j \leq T}
\Bigg[ z_{i}^{T-j} \Bigg]
= \prod_{1 \leq i < j \leq T} (z_{i}-z_{j}).
\label{eqn:Vandermonde}
\end{equation}
\end{lem}
}}
\vskip 0.5cm

\noindent This formula clarifies that the Schur functions are
{\bf symmetric polynomials} in $\z=(z_{1}, \cdots, z_{T})$.

For the second formula, we define the polynomials
$e_{j}(z_{1}, \cdots, z_{T})$'s as the 
coefficients in the expansion
\begin{equation}
\prod_{i=1}^{T} (1 + z_{i} \xi) 
= \sum_{j=0}^{T} e_{j}(z_{1}, \cdots, z_{T}) \xi^{j}.
\label{eqn:defe}
\end{equation}
Then
\begin{eqnarray}
e_{j}(z_{1}, \cdots, z_{T}) &=&
\mbox{ sum of all monomials in the form
$z_{i_{1}} z_{i_{2}} \cdots z_{i_{j}}$} \nonumber\\
&& \mbox{for all strictly increasing sequences 
$1 \leq i_{1} < i_{2} < \cdots < i_{j} \leq T$}.
\nonumber
\end{eqnarray}
$e_{j}(z_{1}, \cdots, z_{T})$'s are also symmetric
polynomials in $z_{1}, \cdots, z_{T}$ and called
the {\bf $j$-th elementary symmetric polynomials}.
\vskip 0.5cm
\fbox{\parbox{16.5cm}{
\begin{lem}
\label{thm:JT2}
Assume that the conjugate of $\vlambda$ is given by
$\widetilde{\vlambda}=(\widetilde{\lambda}_{1}, \cdots,
\widetilde{\lambda}_{N})$ with length $N$. Then
$$
s_{\lambda}(z_{1}, \cdots, z_{T})
= \det_{1 \leq i, j \leq N}
\Bigg[ e_{\widetilde{\lambda}_{j}+(i-j)}
(z_{1}, \cdots, z_{T}) \Bigg].
$$
\end{lem}
}}
\vskip 0.5cm
\noindent More details for YD, SSYT and 
symmetric polynomials, see {\it e.g.}
Fulton (1997) \cite{Ful97}.

Now we go back to the vicious walker model.
We notice a simple relation between the
partition ${\bf L}=(L_{1}, \cdots, L_{N})$
and the final positions of the $N$ vicious walkers
at time $T$,
$
\y=(y_{1}, \cdots, y_{N}) \equiv
(S_{1}(T), \cdots, S_{N}(T)),
$
given by
$
   y_{i}=T-2 L_{i}+2(i-1), \quad 1 \leq i \leq N,
$
on the initial condition (\ref{eqn:initial}).
Then we have established the following
relation between the vicious walks and
YD/SSYT/Schur functions.
\vskip 0.5cm
\fbox{\parbox{16.5cm}{
Set the initial positions as $\x_{0}=(0,2, \cdots, 2(N-1))$.
For $\y$ such that $\y \in \Z_{<}^{N}$, if
$T \in 2 \N$,
$\y+1 \equiv (y_{1}+1, \cdots, y_{N}+1) 
\in \Z_{<}^{N}$, if $T+1 \in 2 \N$,
let 
${\bf L}=(L_{1}, \cdots, L_{N})$ with
$L_{i}=(T+x_{0i}-y_{i})/2, \, 1 \leq i \leq N,$
and $\vlambda=\widetilde{\bf L}$. Then
\begin{eqnarray}
&& \qquad \bullet \quad
\mbox{positions $\y$ of vicious walkers} \nonumber\\
&& \qquad \qquad 
\mbox{at the finial time $T$} \qquad \qquad \quad
\Longleftrightarrow \quad
\mbox{YD with the shape $\vlambda$} \nonumber\\
&& \qquad \bullet \quad
\mbox{ a realization of vicious walk} \nonumber\\
&& \qquad \qquad
\mbox{from $\S(0)=\x_{0}$ to $\S(T)=\y$} \quad
\Longleftrightarrow \quad
\mbox{an SSYT ${\cal T}$ with the shape $\vlambda$} \nonumber\\
&& \qquad \bullet \quad
\mbox{a set of all vicious walks} \nonumber\\
&& \qquad \qquad \mbox{from $\S(0)=\x_{0}$ to $\S(T)=\y$} \quad
\Longleftrightarrow \quad
\mbox{a Schur function} \,
s_{\lambda}(z_{1}, \cdots, z_{T}).
\nonumber
\end{eqnarray}
}}
\vskip 0.5cm
For $\y=(y_{1}, y_{2}, \cdots, y_{N})$ such that
$\y \in \Z_{<}^{N}$, if $T \in 2 \N$,
$\y+1 \in \Z_{<}^{N}$, if $T+1 \in 2 \N$, and
$\x_{0}=(0,2, \cdots, 2(N-1)),$
define the number
\begin{eqnarray}
M_{N,T}(\y) &=& \mbox{total number of distinct realizations of
vicious walk of $N$ walkers} \nonumber\\
&& \mbox{from the positions 
$\S(0)=\x_{0}$
to $\S(T)=\y$}. 
\nonumber
\end{eqnarray}
The above relations prove the following
identity.
\vskip 0.5cm
\fbox{\parbox{16.5cm}{
Assume that $L_{i}=\displaystyle{\frac{T+2(i-1)-y_{i}}{2}},
{\bf L}=(L_{1}, \cdots, L_{N})$, and
$\vlambda=\widetilde{\bf L}$. Then
$$
M_{N,T}(\y)= s_{\lambda}(z_{1}, \cdots, z_{T}) 
\Bigg|_{z_{1}=z_{2}=\cdots=z_{T}=1}.
$$
}}
\vskip 0.5cm
Lemma \ref{thm:JT1} gives the following estimate
for $M_{N,T}(\y)$, via appropriate
{\bf $q$-factorization} and the formula
for the Vandermonde determinant (\ref{eqn:Vandermonde}).
\vskip 0.5cm
\fbox{\parbox{16.5cm}{
\begin{prop}
\label{thm:prop3}
\begin{eqnarray}
M_{N,T}(\y) &=& \lim_{q \to 1} s_{\lambda}
(1,q, q^2, \cdots, q^{T-1}) \nonumber\\
&=& \lim_{q \to 1} q^{\sum_{k=1}^{T} (k-1) \lambda_{k}}
\prod_{1 \leq i < j \leq T}
\frac{q^{\lambda_{i}-\lambda_{j}+j-i}-1}
{q^{j-i}-1} 
= \prod_{1 \leq i < j \leq T}
\frac{\lambda_{i}-\lambda_{j}+j-i}{j-i}.
\nonumber
\end{eqnarray}
\end{prop}
}}
\vskip 0.3cm
On the other hand, 
by definition (\ref{eqn:defe}), it is easy to see that
$$
e_{j}(z_{1}, \cdots, z_{T}) \Bigg|_{z_{1}=\cdots=z_{T}=z}
= {T \choose j} z^{T}
\equiv \frac{T!}{j ! (T-j)!} z^{T}.
$$
Then Lemma \ref{thm:JT2} gives the following.
\vskip 0.5cm
\fbox{\parbox{16.5cm}{
\begin{prop}
\label{thm:prop4}
\begin{eqnarray}
M_{N,T}(\y) &=& \det_{1 \leq i, j \leq N}
\left[ 
{T \choose \widetilde{\lambda}_{j}+i-j} \right]
= \det_{1 \leq i, j \leq N}
\left[ 
{T \choose L_{j}+i-j} \right]
\nonumber\\
&=& \det_{1 \leq i, j \leq N}
\left[ 
{T \choose \{T+2(i-1)-y_{j}\}/2 } \right]. \nonumber
\end{eqnarray}
\end{prop}
}}

\setcounter{equation}{0}
\section{Determinantal Formula for Nonintersecting Paths}
Since we have assumed the initial positions as
(\ref{eqn:initial}), we can see that
the $(i, j)$-element of the matrix
in the determinant in Proposition \ref{thm:prop4}
is
\begin{eqnarray}
{T \choose \{T+2(i-1)-y_{j}\}/2 }
&=& \# \Big\{ \mbox{{\bf lattice path}
from $2(i-1)$ at time 0 to $y_{j}$ at time $T$} 
\Big\}
\nonumber\\
&=& \sum_{\mbox{\footnotesize 
all lattice paths: $(2(i-1),0) \leadsto (y_{j},T)$}} 
w(\mbox{path}) \Bigg|_{w(\mbox{\scriptsize path})=1}.
\nonumber
\end{eqnarray}
If we define an appropriate weight function $w(\mbox{path})$
on single lattice paths, the summation of $w(\mbox{path})$
will give the {\bf Green function} of single lattice paths,
$$
G \Big((x,0), (y,T) \Big)=\sum_{\mbox{\footnotesize 
all lattice paths: $(x,0) \leadsto (y,T)$}} w(\mbox{path}).
$$
Proposition \ref{thm:prop4} can be regarded as a special
case of the {\bf Karlin-McGregor formula} 
in the probability theory \cite{KM59a,KM59b},
and the {\bf Lindstr\"om-Gessel-Viennot formula} in the 
enumerative combinatorics (see \cite{Ste90,KGV00} and references therein). 
In order to explain this fact, here we introduce some
definitions and notations for describing 
{\bf lattice paths}.

Let $V=\{ \mbox{vertex}\}, E=\{\mbox{directed edge}\}$, 
$D=(V,E)=$ an {\bf acyclic} directed graph, 
where acyclic means that any cycles of directed edges are forbidden.
For $u, v \in V$, 
\begin{eqnarray} 
\mbox{a {\bf lattice path}} \ u \to v &=& 
\mbox{a sequence of directed edges from $u$ to $v$}, \nonumber\\
{\cal P}(u,v) &=& \ \mbox{the set of all lattice paths from
$u$ to $v$}. \nonumber
\end{eqnarray}
A weight function 
$w: E \to \Z[[x_{e}: e \in E]]$ 
is introduced,
where $\Z[[x_{e}: e \in E]]$ denotes 
a ring of formal power series
of $\{x_{e} : e \in E\}$ and the
{\bf weight on a lattice path} $P$ is defined by
$\displaystyle{
w(P) \equiv \prod_{e \in P}w(e). }$
Then the
{\bf Green function of lattice paths} 
from $u$ to $v$ is defined by
$$
   G(u, v)= \sum_{P: P \in {\cal P}(u,v)} w(P).
$$
Let $I =\{u_{1}, u_{2}, \cdots, u_{N}\}, 
J=\{v_{1}, v_{2}, \cdots, v_{N}\}$ with 
$ u_{i}, v_{i} \in V, i=1,2, \cdots, N.$
The sets $I, J$ are ordered;
$u_{1} < u_{2} < \cdots < u_{N}, v_{1} < v_{2} < \cdots < v_{N}$.
Then we consider a set of 
{\bf $N$-tuples of lattice paths}
$$
 {\cal P}(I,J)=\Big\{\P=(P_{1}, \cdots, P_{N}):
  P_{i} \in {\cal P}(u_{i}, v_{i}), i=1,2, \cdots, N \Big\}.
$$
The weight for each $N$-tuple of lattice paths is given by
$$
  w(\P)=\prod_{i=1}^{N} w(P_{i}) 
  = \prod_{i=1}^{N} \prod_{e \in P_{i}} w(e).
$$
We say `lattice paths $P$ and $Q$ {\bf intersect}', if 
$P$ and $Q$ share at least one common vertex.
Then, for ordered sets of vertices $I$ and $J$, we say 
`$I$ is $D$-{\bf compatible} with $J$'
in the case that,
whenever $u_{i} < u_{j}$ in $I$ and $v_{i} < v_{j}$ in $J$, 
every lattice path $P \in {\cal P}(u_{i}, v_{j})$ intersects
every lattice path $Q \in {\cal P}(u_{j}, v_{i})$.
A set of {\bf $N$-tuples of nonintersecting 
lattice paths} is denoted by
$$
{\cal P}_{0}(I, J) = \Big\{\P \in {\cal P}(I,J) :
\mbox{any lattice 
paths in $\P$ do not intersect with others} \Big\},
$$
and the {\bf Green function of $N$-tuples of nonintersecting 
lattice paths}
is defined by
$$
G_{\rm nonint}(I,J)=\sum_{\P \in {\cal P}_{0}(I,J)}
w(\P).
$$

\vskip 0.5cm
\fbox{\parbox{16.5cm}{
\begin{thm}
\label{thm:Gessel-Viennot}
Let $I=(u_{1}, \cdots, u_{N})$ and
$J=(v_{1}, \cdots, v_{N})$ be two ordered sets of 
vertices in an acyclic graph $D$.
If $I$ is $D$-compatible with $J$, then
the Green function of the nonintersecting $N$-tuples
of lattice paths is given by
$$
G_{\rm nonint}(I,J)
= \det_{1 \leq i, j \leq N}
\Bigg[ G(u_{i}, v_{j}) \Bigg],
$$
where $G(u, v)$ denotes the Green function of
single lattice paths from $u$ to $v$ on $D$.
\end{thm}
}}
\vskip 0.5cm

\noindent The proof of Theorem \ref{thm:Gessel-Viennot}
is given in {\sf Appendix A}
following Stembridge (1990) \cite{Ste90}.
For our vicious walker model,
consider the directed graph $D=(V, E)$, where
$$
V=\Big\{(x,t) \in \Z^{2}: x+t = {\rm even}, 
\ t=0,1,2, \cdots, T \Big\},
$$
and all edges connecting the nearest-neighbor pairs
of vertices in $V$ are oriented to the positive direction of $t$ axis.
Set the weight function 
\begin{eqnarray}
w(e)=
\left\{
   \begin{array}{ll}
      1 & \mbox{for} \ e=\langle (x, t-1) \to (x-1,t) \rangle \\
      1 & \mbox{for} \ e=\langle (x, t-1) \to (x+1,t) \rangle \\
      0 & \mbox{otherwise.} \\
   \end{array}\right.
\nonumber
\end{eqnarray}
Set $T \in \N$ and $u_{i}=(x_{i}, 0), v_{i}=(y_{i}, T)
\in V, i=1,2,\cdots, N$. Then the Green function of single 
lattice paths
from $u_{i}$ to $v_{i}$ is 
$$
G(u_{i}, v_{i}) = \Big| {\cal P}(u_{i} \to v_{i}) \Big|
= { T \choose (T+x_{i}-y_{i})/2}.
$$
Theorem \ref{thm:Gessel-Viennot} then gives
the Green function of the $N$-tuples of
nonintersecting lattice paths from
$I=\Big\{(x_{i}, 0) \Big\}_{i=1}^{N}$ to
$J=\Big\{(y_{i}, T) \Big\}_{i=1}^{N}$ as
$\displaystyle{
\det_{1 \leq i, j \leq N} \left[
{T \choose (T+x_{i}-y_{j})/2} \right].
}$

For $\x \in \Z_{<}^{N}$, and $\y$ such that
$\y \in \Z_{<}^{N}$, if $T \in 2 \N$,
$\y+1 \in \Z_{<}^{N}$, if $T+1 \in 2 \N$,
define 
\begin{eqnarray}
M_{N}(T,\y|\x) &=& \mbox{total number of distinct realizations of
vicious walk of $N$ walkers} \nonumber\\
&& \mbox{from the positions $\x=(x_{1},\cdots, x_{N})$
to $\y=(y_{1}, \cdots, y_{N})$ during time $T$}. 
\nonumber
\end{eqnarray}
Proposition \ref{thm:prop4} is now generalized
as follows.

\vskip 0.3cm
\fbox{\parbox{16.5cm}{
\begin{prop}
\label{thm:prop6} \qquad
$\displaystyle{
M_{N}(T,\y|\x)
=\det_{1 \leq i, j \leq N}
\left[ 
{T \choose (T+x_{i}-y_{j})/2 } \right]. 
}$
\end{prop}
}}
\setcounter{equation}{0}
\section{Diffusion Scaling Limit}
Recall that $(\{\S(t)\}_{t=0,1,2,...,T}, {\sf Q}^{\x}_{T})$ 
denotes the vicious walk 
with the noncolliding condition up to time
$T>0$ starting from the positions $\x \in \Z_{<}^{N}$.
For $L \geq 1$, 
we consider probability measures $\mu_{L,T}^{\x}$
on the space of continuous paths $C( [0,T]\to \R^N)$ 
defined by
$$
\mu_{L,T}^{\x}(\cdot) 
= {\sf Q}_{L^2T}^{\x} 
\left(\frac{1}{{L}}\S(L^2t) \in \cdot \right),
$$
where $\S(t), t\ge 0$, is now considered to be
the interpolation of the $N$-dimensional random walk
$\S(t), t=0,1,2,\dots$.
We study the limit of the probability measure
$\mu_{L,T}^{\x}$, $L\to \infty$.

We put 
$$
\RV=\Big\{\x\in \R^N : x_1 < x_2 < \cdots < x_N \Big\},
$$
which can be called the {\bf Weyl chamber of type} $A_{N-1}$ 
(see, for example, \cite{FH91}).
By virtue of the Karlin-McGregor formula \cite{KM59a,KM59b}, 
the transition density function  $f_N (t, \y|\x)$ of
the {\bf absorbing Brownian motion} in $\RV$
and the probability $\cN (t, \x)$ that the Brownian motion 
starting from $\x\in\RV$ does not
hit the boundary of $\RV$ up to time $t>0$
are given by
$$
f_{N}(t, \y|\x)= 
\det_{1 \leq i, j \leq N}
\left[ \frac{1}{\sqrt{2 \pi t}} \ e^{-(x_{j}-y_{i})^2/2t} \right],
\: \x,\y \in \RV,
$$
and
$\displaystyle{
\cN (t, \x ) = 
\int_{\RV} d\y f_N (t, \y |\x ),}
$
respectively. We put 
$\displaystyle{
h_N(\x) = \prod_{1\le i <j \le N}( x_j - x_i),}
$
and let $\0=(0,0,\cdots, 0) \in \R^{N}$, 
which describes the state that
all $N$ particles are at the origin $0$.

\begin{figure}[bhtp]
\includegraphics[width=0.3\linewidth]{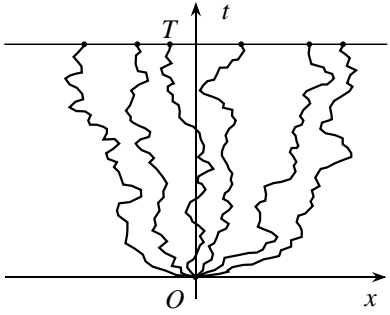}
\caption{ The process $\X(t), t \in [0,T]$, starting from
$\0$.
\label{fig:Fig4}}
\end{figure}

\vskip 0.5cm
\fbox{\parbox{16.5cm}{
\begin{thm}
\label{thm:main}
{\rm (i)}
For any fixed $\x \in \ZV$ and $T>0$, as $L\to \infty$,
$\mu_{L,T}^{\x}(\cdot)$ converges weakly to the law of 
the {\bf temporally inhomogeneous diffusion process}
${\bf X}(t)=(X_1(t),X_2(t),\dots,X_N(t)), t\in [0,T]$,
with transition probability density $g_{N,T}(s,\x; t,\y)$;
\begin{eqnarray}
&& g_{N,T}(0, {\bf 0}; t, \y)
=c_{N}T^{N(N-1)/4}t^{-N^2/2}  
\exp\left\{ -\frac{|\y|^2}{2t} \right\}
h_N(\y)\cN (T-t,\y), \nonumber\\
&&
g_{N,T}(s,\x; t, \y )
= \frac{f_{N}(t-s, \y|\x)\cN (T-t,\y)}{\cN (T-s,\x)},
\nonumber
\end{eqnarray}
for $0 < s < t \le T,\; \x, \y \in \RV,$
where 
$c_N = 2^{-N/2}/
\prod_{i=1}^N \Gamma(i/2)$
with the gamma function $\Gamma$.

\noindent
{\rm (ii)}
The diffusion process ${\bf X}(t)$
solves the following equation:
$$
dX_i (t) = dB_i(t) + b_i^T(t, {\bf X}(t)) dt,
\quad t\in [0,T],\quad i=1,2,\dots, N,
$$
where $B_i(t)$, $i=1,2,\dots,N,$ 
are independent one-dimensional Brownian motions and
$$
b_i^T (t,\x) = \frac{\partial}{\partial x_i}\ln \cN (T-t, \x),
\quad i=1,2,\dots,N.
$$
\end{thm}
}}
\vskip 0.3cm
\noindent Figure \ref{fig:Fig4} illustrates the process
$\X(t), t \in [0,T] $, when all $N$ particles start from the
origin; $\X(0)=\0$.

\vskip 0.3cm
Although $\mu_{L,T}^{\x}(\cdot)$ is the probability measure
defined on $C([0,T]\to \R^N)$,
it can be regarded as that on $C([0,\infty)\to \R^N)$
concentrated on the set
$
\{ w \in C([0,\infty)\to \R^N) : w(t)= w(T),\: t \ge T\}.
$
Next we consider the case that $T=T(L)$ goes to infinity 
as $L \to \infty$.

\vskip 0.3cm
\fbox{\parbox{16.5cm}{
\begin{cor}
\label{thm:VI}
{\rm (i)}
Let $T(L)$ be an increasing function of $L$ with
$T(L) \to \infty$ as $L\to\infty$.
For any fixed $\x \in \ZV$, as $L\to \infty$,
$\mu_{L,T(L)}^{\x}(\cdot)$ converges weakly to the law of 
the {\bf temporally homogeneous diffusion process}
${\bf Y}(t)=(Y_1(t),Y_2(t),\dots,Y_N(t)), t\in [0,\infty)$,
with transition probability density $p_N(s,\x;t,\y)$;
\begin{eqnarray}
&& p_{N}(0, {\bf 0}; t, \y)
= c^{\prime}_{N}t^{-N^2/2} 
\exp\left\{ -\frac{|\y|^2}{2t} \right\} h_N(\y)^2, \nonumber\\
&&
\label{eqn:pN0}
p_N(s,\x; t, \y )
= \frac{1}{h_N(\x)}f_{N}(t-s, \y|\x)h_N(\y),
\end{eqnarray}
for $0 < s < t < \infty,\; \x, \y \in \RV,$
where 
$c^{\prime}_{N}=(2\pi)^{-N/2}
/\prod_{i=1}^{N} \Gamma(i)$.

\noindent
{\rm (ii)}
The diffusion process ${\bf Y}(t)$ solves 
the equations of {\bf Dyson's Brownian motion model}
with the parameter $\beta=2$,
$$
dY_i (t) = dB_i(t) 
+ \sum_{1\le j \le N,j\not=i} \frac{1}{Y_i(t)-Y_j(t)}dt,
\quad t\in [0,\infty), \quad i=1,2,\dots, N.
$$
\end{cor}
}}
\vskip 0.3cm

Here we only give the proof of a key lemma used to prove
Theorem \ref{thm:main}, in order to demonstrate that
the Schur function plays an important role.
See Katori and Tanemura \cite{KT03}
for the complete proofs of Theorem \ref{thm:main}
and Corollary \ref{thm:VI}.

For $L >0$ we introduce the following functions:
$$
\phi_L(x)= 2\left[\frac{{L}x}{2}\right], \ x \in \R,
\hbox{ and }
\phi_L(\x)= \Big(\phi_L(x_1), \phi_L(x_2),
\dots, \phi_L(x_N) \Big),
\ \x \in \R^N,
$$
where $[a]$ denotes the largest integer not greater than $a$.
Let
$
V_{N}(T,\y|\x)=2^{-NT} M_{N}(T,\y|\x),
$
where $M_{N}(T,\y|\x)$ is given by Proposition \ref{thm:prop6}.

\vskip 0.5cm
\fbox{\parbox{16.5cm}{
\begin{lem}
\label{thm:lemma9}
For $t>0$, $\x \in \ZV$ and $\y\in \RV$.
\begin{eqnarray}
&& \left(\frac{{L}}{2}\right)^{N} 
V_N \Big(\phi_{L^2}(t), \phi_L(\y) \Big| \x \Big) 
\nonumber\\
&& \qquad \qquad
= c_N' t^{-N^2/2}
h_N\left( \frac{\x}{L}\right)
\exp \left\{ -\frac{|\y|^2}{2t}\right\} h_N(\y)
 \left( 1+ {\cal O}\left(\frac{|\y|}{L}\right) \right),
\nonumber
\end{eqnarray}
as $L \to \infty$.
\par\noindent
\end{lem}
}}
\vskip 0.3cm
\noindent{\sf Proof.} 
It will be enough to consider the case that
$\x = 2 \u =(2u_{1}, \cdots, 2u_{N}) \in \Z_{<}^{N}, 
\y=2 \v =(2 v_{1}, \cdots, 2 v_{N}) \in \Z_{<}^{N}$,
and
$\phi_{L^2}(t) =2\ell, \ell \in \Z_{+}$,
where $\Z_{+}=\{1,2,3, \cdots\}$.
Then 
$$
M_N \Big(\phi_{L^2}(t) , \phi_L(\y) \Big| \x \Big) 
=
M_N ( 2\ell , 2\v |2\u ) 
= \det_{1 \leq i, j \leq N}
\left[ 
{2\ell \choose \ell +u_j-v_i}
\right],
$$
and
\begin{eqnarray}
{2\ell \choose \ell +u_j-v_i}
&=&
\frac{(2\ell)!}{(\ell+u_j-v_i)! (\ell-u_j+v_i)!}
\nonumber
\\
&=&
\frac{(2\ell)!}{(\ell-v_i)! (\ell+v_i)!}A_{ij}(\ell,\v,\u),
\nonumber
\end{eqnarray}
with
$$
A_{ij}(\ell,\v,\u)=
\frac{(\ell+v_i-u_j+1)_{u_j}}{(\ell-v_i+1)_{u_j}},
$$
where we have used the Pochhammer symbol;
$(a)_0 \equiv 1$,
$(a)_i= a(a+1) \cdots (a+i-1)$, $i \ge 1$.
Then
\begin{equation}
M_N \Big(\phi_{L^2}(t), \phi_L(\y) \Big| \x \Big) 
= \prod_{i=1}^N
\frac{(2\ell)!}{(\ell-v_i)! (\ell+v_i)!}
\det_{1 \leq i, j \leq N}
\Bigg[ 
A_{ij}(\ell,\v,\u)
\Bigg].
\label{eq:Nn}
\end{equation}
The leading term of 
$\displaystyle{
\det_{1 \leq i, j \leq N} \Bigg[ A_{ij}(\ell,\v,\u) \Bigg]} $ 
in $L\to\infty$ is
\begin{eqnarray}
D_1(\v,\u) 
&=& 
\det_{1 \leq i, j \leq N}
\left[ \left(\frac{\ell +v_i}{\ell -v_i}\right)^{u_j} \right]
\nonumber
\\
&=&
(-1)^{N(N-1)/2}
\det_{1 \leq i, j \leq N}
\left[ \left(\frac{\ell +v_i}{\ell -v_i}\right)^{u_{N-j+1}} \right].
\nonumber
\end{eqnarray}
Let $\vxi(\u)=(\xi_{1}(\u), \dots, \xi_{N}(\u))$ 
be a partition specified by the starting point
$2 \u$ defined by
$$
 \xi_{j}(\u)=u_{N-j+1}-(N-j), \ j =1,2,\dots, N.
$$
We have
\begin{eqnarray}
D_1(\v,\u) 
&&= 
(-1)^{N(N-1)/2}
\det_{1 \leq i, j \leq N}
\left[ \left(\frac{\ell +v_i}{\ell -v_i}\right)^{N-j} \right]
s_{\xi(\u)}
\left(\frac{\ell+v_1}{\ell-v_1},\dots,\frac{\ell+v_N}{\ell-v_N} \right)
\nonumber
\\
&&=
(-1)^{N(N-1)/2}
\prod_{1\le i<j\le N}
\left( \frac{\ell+v_i}{\ell- v_i}- \frac{\ell+v_j}{\ell-v_j} \right)
s_{\xi(\u)}
\left(\frac{\ell+v_1}{\ell-v_1},\dots,\frac{\ell+v_N}{\ell-v_N} \right)
\nonumber
\\
&&=
\prod_{1\le i<j\le N}
\frac{2\ell(v_j-v_i)}{(\ell-v_i)(\ell-v_j)}
s_{\xi(\u)}
\left(\frac{\ell+v_1}{\ell-v_1},\dots,\frac{\ell+v_N}{\ell-v_N} \right),
\nonumber
\end{eqnarray}
where we have used Lemma \ref{thm:JT1} 
for the Schur function associated to the partition $\vxi(\u)$.
Proposition \ref{thm:prop3} gives
$$
s_{\xi(\u)}(1,1, \dots, 1)=
\prod_{1 \leq i < j \leq N}
\frac{\xi_{i}(\u)-\xi_{j}(\u)+j-i}{j-i}.
$$
Therefore the leading term of $D_1(\v,\u)$ in $L\to\infty$ is
\begin{eqnarray}
D_2(\v,\u) 
&=&
\prod_{1\le i<j\le N}
\frac{2(v_j-v_i)}{\ell} \times
s_{\xi(\u)}(1,1,\dots,1)
\nonumber
\\
&=&
\ell^{- N(N-1)/2}2^{N(N-1)/2}h_N(\v)h_N(\u)
\prod_{1\le i<j\le N}\frac{1}{j-i}
\nonumber
\\
&=&
h_N\left(\frac{\v}{\ell}\right)h_N(2\u)
\prod_{i=1}^N \frac{1}{\Gamma(i)}.
\label{eq:D2}
\end{eqnarray}
On the other hand, 
by Stirling's formula we see that
\begin{equation} 
\prod_{i=1}^N
\frac{(2\ell)!}{(\ell-v_i)! (\ell+v_i)!}
= (\ell \pi)^{-N/2} 2^{2N\ell}
\prod_{i=1}^N 
\left( 1- \frac{v_i^2}{\ell^2}\right)^{-\ell-1/2}
\left( \frac{1-v_{i}/\ell}{1+v_{i}/\ell}\right)^{v_i}
\left(1+{\cal O}\left(\frac{1}{\ell}\right)\right).
\label{eq:St}
\end{equation}
From (\ref{eq:Nn}), (\ref{eq:D2}) and (\ref{eq:St})
\begin{eqnarray}
&&V_N \Big(\phi_{L^2}(t), \phi_L(\y) \Big| \x \Big)
= 2^{-2N\ell} 
M_N \Big(\phi_{L^2}(t), \phi_L(\y) \Big| \x \Big)
\nonumber
\\
&& =
c_N' \left(\frac{2}{\ell}\right)^{N/2}
h_N\left(\frac{\v}{\ell}\right)
h_N\left( 2 \u \right) 
\exp \left\{ - \frac{|\v|^2}{\ell} \right\}
\left(1+{\cal O}\left(\frac{|\v|}{\ell}\right)\right)
\nonumber
\\
&&= 
c_N' \left(\frac{2}{L}\right)^N t^{-N^2/2}
h_N\left(\frac{\x}{L}\right)
\exp \left\{- \frac{|\y|^2}{2t}\right\}
h_N(\y)\left(1+{\cal O}\left(\frac{|\y|}{L}\right)\right).
\nonumber
\end{eqnarray}
Then we obtain Lemma \ref{thm:lemma9}. \qed
\vskip 0.5cm

Corollary \ref{thm:VI} is obtained from Theorem \ref{thm:main}
by the following evaluation of asymptotic \cite{KT03}.
Let $t > 0$ and $\x \in \RV$, then
$$
{\cal N}_{N}(t, \x)
=\frac{1}{\overline{c}_N}
h_N \left( \frac{\x}{\sqrt{t}}\right)
\left( 1+ {\cal O}\left(\frac{|\x|}{\sqrt{t}}\right) \right),
\quad \mbox{in the limit} \quad \frac{|\x|}{\sqrt{t}} \to 0,
$$
where $\overline{c}_{N}=\pi^{N/2} \prod_{i=1}^{N}
\Big\{\Gamma(i)/\Gamma(i/2) \Big\}$.

\setcounter{equation}{0}
\section{Eigenvalue Process of Hermitian Matrix-valued Processes}

We consider complex-valued processes 
$\xi_{ij}(t) \in \C, 1 \leq i, j \leq N, t \in [0, \infty),$
with the condition $\xi_{ji}(t)^{*}=\xi_{ij}(t)$,
and introduce {\bf Hermitian matrix-valued processes} 
$\Xi(t)= \Big(\xi_{ij}(t) \Big)_{1 \leq i, j \leq N}$.
We denote by 
$U(t)=\Big(u_{ij}(t) \Big)_{1 \leq i, j \leq N}$ 
the family of unitary matrices which diagonalize $\Xi(t)$ so that
\begin{equation}
U(t)^{\dagger} \Xi (t) U(t)=\Lambda(t)
={\rm diag} \Big\{\lambda_{1}(t), \lambda_{2}(t), \cdots,
\lambda_{N}(t) \Big\},
\label{eqn:diag}
\end{equation}
where $\Big\{\lambda_{i}(t) \Big\}_{i=1}^{N}$ 
are eigenvalues of $\Xi(t)$ and we assume their increasing order
$$
\lambda_{1}(t) \leq \lambda_{2}(t) \leq \cdots 
\leq \lambda_{N}(t).
$$
Define $\Gamma_{ij}(t), 1 \leq i, j \leq N$, by
\begin{equation}
\Gamma_{ij}(t) dt = 
\Big(U(t)^{\dagger} d \Xi(t) U(t) \Big)_{ij}
\Big(U(t)^{\dagger} d \Xi(t) U(t) \Big)_{ji},
\label{eqn:Gamma}
\end{equation}
where $d \Xi(t)=\Big(d \xi_{ij}(t) \Big)_{1 \leq i, j \leq N}$.
The indicator function ${\bf 1}_{\{\omega\}}$ gives
${\bf 1}_{\{\omega\}}=1$ if the condition $\omega$ is satisfied, and
${\bf 1}_{\{\omega\}}=0$ otherwise.

\vskip 0.3cm
\fbox{\parbox{16.5cm}{
\begin{thm}
\label{thm:BruThm}
Assume that $\xi_{ij}(t), 1 \leq i < j \leq N$, are
continuous semimartingales.
The process of eigenvalues 
$\vlambda(t)=(\lambda_{1}(t), \lambda_{2}(t),
\cdots, \lambda_{N}(t))$ satisfies 
the {\bf stochastic differential equations}
\begin{equation}
d \lambda_{i}(t)=dM_{i}(t)+d J_{i}(t), \quad
 t \in [0, \infty), \,
i=1,2, \cdots, N, \,
\label{eqn:SDE1}
\end{equation}
where $M_{i}(t)$ is the martingale with quadratic variation
\begin{equation}
  \langle M_{i} \rangle_{t}=
  \int_{0}^{t} \Gamma_{ii}(s) ds
\label{eqn:SDE2}
\end{equation}
and $J_{i}(t)$ is the process with finite variation given by
\begin{equation}
dJ_{i}(t) = \sum_{j=1}^{N} \frac{1}{\lambda_{i}(t)-\lambda_{j}(t)}
{\bf 1}_{\{\lambda_{i}(t) \not= \lambda_{j}(t)\}} \Gamma_{ij}(t) dt
+ d \Upsilon_{i}(t)
\label{eqn:SDE3}
\end{equation}
where $d \Upsilon_{i}(t)$ is the finite-variation part of
$\Big(U(t)^{\dagger} d \Xi(t) U(t) \Big)_{ii}$.
\end{thm}
}}
\vskip 0.5cm

\noindent This theorem is obtained by simple generalization
of Theorem 1 in Bru \cite{Bru89}.
A key point to derive the theorem is applying
the {\bf It\^o rule} for differentiating the product of
matrix-valued semimartingales:
If $X$ and $Y$ are $N \times N$ matrices
with semimartingale elements,
then
$$
d( X^{\dagger} Y)= (dX)^{\dagger} Y + X^{\dagger} (d Y) +
(dX)^{\dagger} (dY).
$$
We give the proof in {\sf Appendix B}.
(See Remark 3 below.)

Let $B_{ij}(t)$, $\widetilde{B}_{ij}(t)$, 
$1\le i,j \leq N$,
be independent one-dimensional Brownian motions. 
For $1 \leq i, j \leq N$ we set
$$
s_{ij}(t)
=
\left\{
   \begin{array}{ll}
      \displaystyle{
      \frac{1}{\sqrt{2}} B_{ij}(t)}, & 
   \mbox{if} \ i < j, \\
        & \\
        B_{ii}(t), & \mbox{if} \ i=j, \\
        & \\
    \displaystyle{
     \frac{1}{\sqrt{2}} B_{ji}(t)}, &
    \mbox{if} \ i > j, \\
   \end{array}\right. 
   \quad {\rm and } \quad
a_{ij}(t)
=
\left\{
   \begin{array}{ll}
      \displaystyle{
      \frac{1}{\sqrt{2}} \widetilde{B}_{ij}(t),
      } & 
   \mbox{if} \ i < j, \\
        & \\
    0, & \mbox{if} \ i=j,\\
        & \\
     \displaystyle{
      -\frac{1}{\sqrt{2}} \widetilde{B}_{ji}(t),
      } & 
   \mbox{if} \ i > j. \\
   \end{array}\right.
$$
A Hermitian matrix-valued process is defined by
\begin{equation}
\Xi(t) = \Big(\xi_{ij}(t) \Big)_{1 \leq i, j \leq N}
=\Big(s_{ij}(t)+\sqrt{-1}a_{ij}(t) \Big)_{1 \leq i, j \leq N},
\quad t \in [0, \infty).
\label{eqn:Hermite1}
\end{equation}
By definition 
$d\xi_{ij}(t) d\xi_{k \ell}(t)=\delta_{i \ell} \delta_{j k} dt,
\,1 \leq i, j, k, \ell \leq N,$
and thus $\Gamma_{ij}(t)=1$. 
Theorem \ref{thm:BruThm} thus implies that
the eigenvalue process
$\vlambda(t)$ of this matrix-valued process (\ref{eqn:Hermite1})
solves the equations of {\bf Dyson's Brownian motion
model} with the parameter $\beta=2$ \cite{Dys62}
\begin{equation}
d \lambda_{i}(t)= dB_{i}(t)+ \frac{\beta}{2}
\sum_{1 \leq j \leq N, j \not=i}
\frac{1}{\lambda_{i}(t)-\lambda_{j}(t)} dt,
\quad t \in [0, \infty), \, i=1,2, \cdots, N,
\label{eqn:Dyson1}
\end{equation}
where $B_{i}(t), i=1,2, \cdots, N$ are
independent one-dimensional Brownian motions.

\begin{description}
\item[Remark 3.] \,
In general, Eqs.(\ref{eqn:SDE1}) with (\ref{eqn:SDE2})
and (\ref{eqn:SDE3}) for the eigenvalue process
$\vlambda(t)=(\lambda_{1}(t), \cdots, \lambda_{N}(t))$
depend on unitary matrix $U(t)$ through 
$\Gamma_{ij}(t)$ defined by (\ref{eqn:Gamma}).
The equations written in the form,
\begin{equation}
d \lambda_{i}(t) =\sum_{j} \alpha_{ij}(t, \vlambda(t))  dB_{j}(t)
+ \beta_{i}(t, \vlambda(t)) dt, 
\quad t \in [0, \infty), \, i=1,2, \cdots, N,
\label{eqn:SDEgeneral}
\end{equation}
where the coefficients $\alpha_{ij}(t, \vlambda)$
and $\beta_{i}(t, \vlambda)$ are functions not only of
$\vlambda$ but also of other variables,
are generally called
{\bf stochastic differential equations (SDE's)} in
\cite{IW89}
(see Definition 1.1 with Eqs.(1.1), (1.1') 
in Chapter IV `Stochastic Differential Equations' on page 159.)
In the special case, in which
these coefficients are only depending on
$\vlambda(t)$, equations are given in the form
\begin{equation}
d \lambda_{i}(t)=\sum_{j} \sigma_{ij}(\vlambda(t)) d B_{j}(t)
+ b_{i}(\vlambda(t)) dt, 
\quad t \in [0, \infty), \, i=1,2, \cdots, N,
\label{eqn:SDEMarkov}
\end{equation}
and they are said to be of the {\bf Markovian type} 
(see page 172 with Eq. (2.11) in \cite{IW89}).
The condition that the SDE's of eigenvalue process 
are reduced to be of the Markovian type 
may be that the matrix-valued process
$\Xi(t)$ is {\bf unitary invariant in distribution}.
By virtue of properties of Brownian motions,
the Hermitian matrix-valued process $\Xi(t)$ 
defined by (\ref{eqn:Hermite1}) is unitary invariant
in distribution, and thus the obtained
SDE's of Dyson's Brownian motion model 
are of the Markovian type.
\end{description}

\setcounter{equation}{0}
\section{Concluding Remarks}

Corollary \ref{thm:VI}(ii) and Eq.(\ref{eqn:Dyson1}) 
with $\beta=2$
implies that the temporally homogeneous
process $\Y(t)$ obtained as a diffusion scaling
limit of vicious walks and the eigenvalue process $\vlambda(t)$
of the Hermitian matrix-valued process (\ref{eqn:Hermite1}) are
equivalent in distribution.
The formula (\ref{eqn:pN0}) in Corollary \ref{thm:VI}
shows that it is the {\bf $h$-transform in the sense of
Doob} \cite{Doo84} of the {\bf absorbing Brownian motion
in the Weyl chamber} $\R_{<}^{N}$, since
$h_{N}(\x)$ is a strictly positive harmonic function
in $\R_{<}^{N}$ \cite{Gra99}.
An interesting relationship between this temporally
homogeneous process $\Y(t)$ (Dyson's Brownian motion model
with the parameter $\beta=2$)
and the temporally inhomogeneous process
$\X(t)$ given by Theorem \ref{thm:main} was reported
in \cite{KT03,KT03b}.
A systematic study on the relations among various
matrix-valued processes, standard, chiral and
non-standard random matrix theories,
and families of noncolliding diffusion processes 
was reported in \cite{KT04}.

For the noncolliding diffusion processes starting from $\0$,
the {\bf multi-time correlation functions} 
were calculated by
Nagao and the present authors using 
the quaternion determinants of self-dual 
quaternion matrices 
({\it i.e.} {\bf pfaffians}) 
and the scaling limits of the infinite particles
$N \to \infty$ and the infinite time-interval
$T \to \infty$ were investigated \cite{NKT03,KNT04}.
Further study of {\bf infinite systems of
noncolliding diffusion particles} 
will be reported elsewhere \cite{KT05}.

\vskip 1cm
\appendix
\setcounter{equation}{0}
\section{Proof of Theorem \ref{thm:Gessel-Viennot}}

By definition of determinant
\begin{equation}
\det_{1 \leq i, j \leq N} \Big[ G(u_{i}, v_{j}) \Big]
= \sum_{\sigma \in \Sym_{N}} {\rm sgn}(\sigma) G(u_{1}, v_{\sigma(1)})
\cdots G(u_{N}, v_{\sigma(N)}),
\label{eqn:LGVp1}
\end{equation}
where $\Sym_{N}$ is the set of all permutations 
of $\{1,2, \cdots, N\}$.
We may interpret (\ref{eqn:LGVp1}) as a generating function for
$(N+1)$-tuples $(\sigma, P_{1}, \cdots, P_{N})$, where
$\sigma \in \Sym_{N}, P_{i} \in {\cal P}(u_{i}, v_{\sigma(i)}), 
i=1,2, \cdots, N$.

Consider an arbitrary configuration $(\sigma, P_{1}, \cdots, P_{N})$
with at least one pair of intersecting 
lattice paths.
We set the order of all vertices in $V$. Let $v$ denote the last
vertex among all vertices that occur as points of intersection
among the lattice paths. Among the lattice paths 
that pass through $v$, 
assume that $P_{i}$ and $P_{j}$ are the two whose indices $i$
and $j$ are the smallest (see Figure 5).

\begin{figure}[htbp]
\begin{minipage}{0.5\linewidth}
\includegraphics[width=0.9\linewidth]{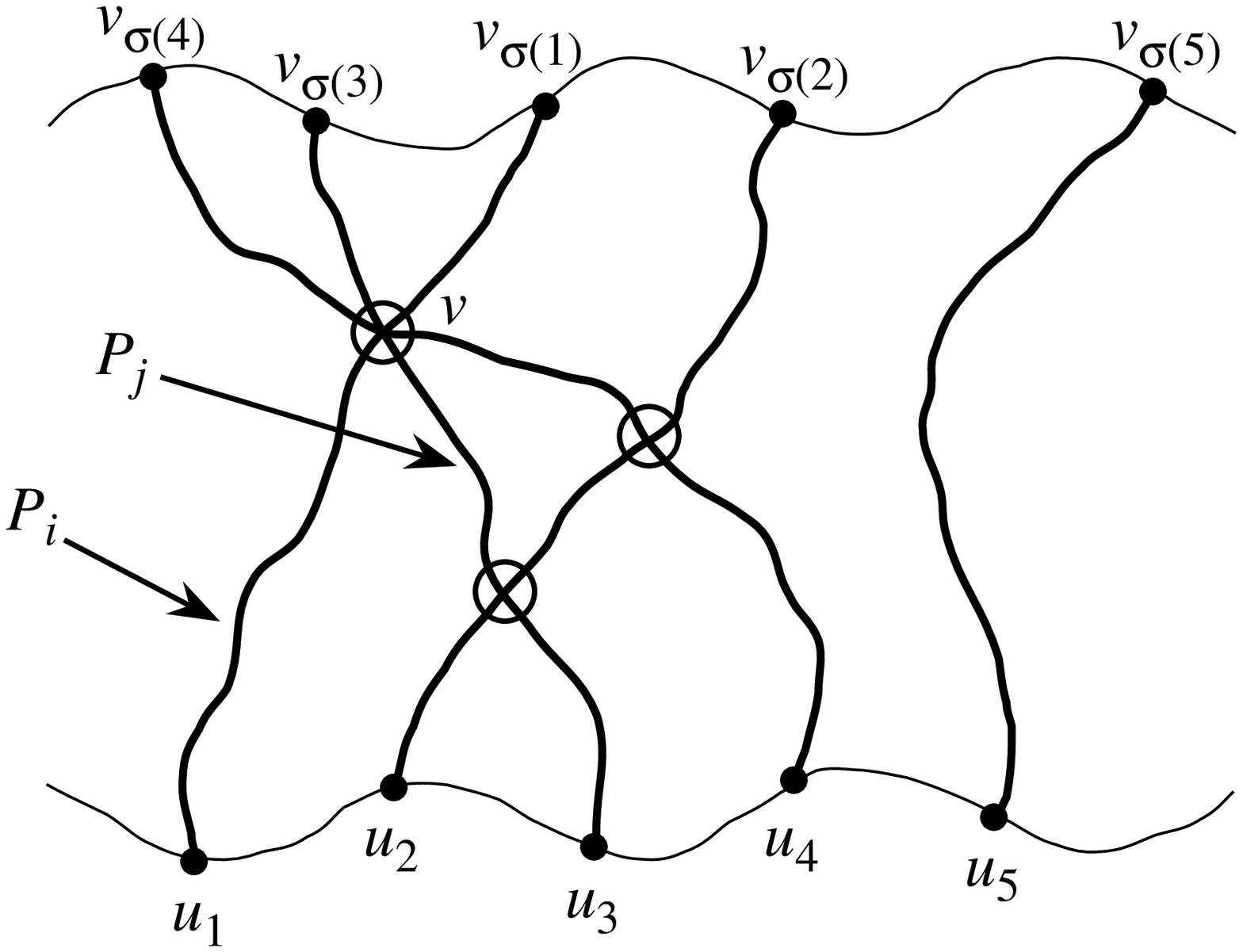}
\label{fig:Fig5}
\vskip 0.5cm
\caption{A configuration $(\sigma, P_{1}, \cdots, P_{N})$.
The last intersecting vertex is denoted by $v$
and lattice paths $P_{i}$ and $P_{j}$ are chosen, 
both of which
pass through $v$.}
\end{minipage}
\begin{minipage}{0.5\linewidth}
\includegraphics[width=0.7\linewidth]{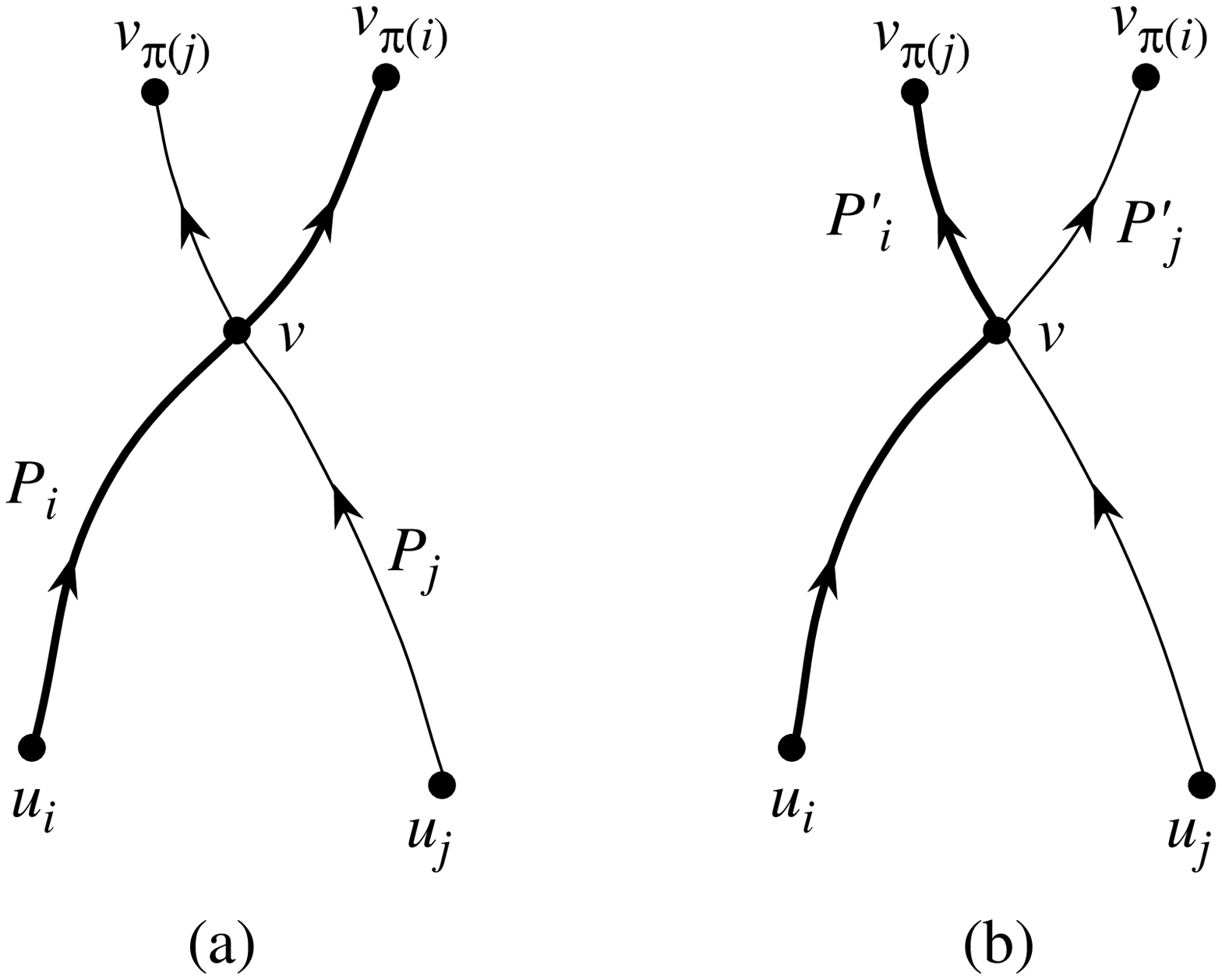}
\caption{(a) Paths $P_{i}$ and $P_{j}$.
(b) Paths $P'_{i}$ and $P'_{j}$.}
\label{fig:Fig6}
\end{minipage}
\end{figure}

Write
$$
    P_{i}=P_{i}(\to v) P_{i}(v \to), \quad \mbox{and} \quad
    P_{j}=P_{j}(\to v) P_{j}(v \to).
$$
For a configuration $(\sigma, P_{1}, \cdots, P_{N})$, define
as shown in Figure 6
\begin{eqnarray}
    P'_{i} &=& P_{i}(\to v) P_{j}(v \to), \nonumber\\
    P'_{j} &=& P_{j}(\to v) P_{i}(v \to), \nonumber\\
    P'_{k} &=& P_{k} \quad \mbox{for} \ k \not=i, j, \nonumber\\
    \sigma' &=& \sigma \circ (i,j), \quad \mbox{where} \
(i, j) \ \mbox{denotes an exchange of $i$ and $j$.} \nonumber
\end{eqnarray}
The operation
$$
   (\sigma, P_{1}, \cdots, P_{N}) 
\mapsto (\sigma', P'_{1}, \cdots, P'_{N})
$$
preserves the set of vertices of intersection, and is
an involution. The weight of lattice paths is 
the same, but the sign is changed.
So any such pair $\Big\{ (\sigma, P_{1}, \cdots, P_{N}),
(\sigma', P'_{1}, \cdots, P'_{N}) \Big\} $ appear in (\ref{eqn:LGVp1})
is canceled out.

The only configurations remain in (\ref{eqn:LGVp1}) are
nonintersecting lattice paths. Since $I$ is assumed to be
$D$-compatible with $J$, for nonintersecting lattice paths
$\sigma={\rm id}$, {\it i.e.},
${\rm sgn}(\sigma)={\rm sgn}({\rm id})=1$.
\hfill \qed

\vskip 0.5cm
\setcounter{equation}{0}
\section{Proof of Theorem \ref{thm:BruThm}}

We consider a matrix-valued process 
$A(t)=(\alpha_{ij}(t))_{1 \leq i, j \leq N}$ defined by
$$
d A(t) =U(t)^{\dagger}dU(t)+ \frac{1}{2} dU(t)^{\dagger} dU(t),
\quad t \in [0, \infty) 
$$
with $A(0)=0$. 
Since $U(t)^{\dagger}U(t)=I_{N}$ for all $t$,
where $I_{N}$ denotes the $N \times N$ unit matrix,
$$
0 = d(U(t)^{\dagger} U(t))=
dU(t)^{\dagger} U(t)+U(t)^{\dagger} dU(t)+dU(t)^{\dagger} dU(t).
$$
Then
\begin{eqnarray}
dA(t)^{\dagger} &=& dU(t)^{\dagger} U(t)+ \frac{1}{2}
dU(t)^{\dagger} dU(t) \nonumber\\
&=& -U(t)^{\dagger} dU(t) - \frac{1}{2} dU(t)^{\dagger} dU(t)
= -dA(t),
\nonumber
\end{eqnarray}
that is, $dA(t)$ is anti-Hermitian. We also see that
\begin{eqnarray}
-dA(t) dA(t) &=& dA(t)^{\dagger} dA(t) \nonumber\\
&=& \left(U(t)^{\dagger} dU(t)+\frac{1}{2} dU(t)^{\dagger}
dU(t) \right)^{\dagger}
\left( U(t)^{\dagger} dU(t)+\frac{1}{2} dU(t)^{\dagger}
dU(t) \right) \nonumber\\
&=& dU(t)^{\dagger} U(t) U(t)^{\dagger} dU(t) 
= dU(t)^{\dagger} dU(t).
\label{eqn:antidA}
\end{eqnarray}
This implies 
\begin{equation}
dU(t)=U(t)\left( dA(t) + \frac{1}{2} dA(t) dA(t) \right).
\label{eqn:dU}
\end{equation}
By (\ref{eqn:diag})
\begin{eqnarray}
d \Lambda(t) &=& dU(t)^{\dagger} \Xi(t) U(t) 
+U(t)^{\dagger} d\Xi(t) U(t) +U(t)^{\dagger} \Xi(t) dU(t) \nonumber\\
&& +dU(t)^{\dagger} d\Xi(t) U(t)+dU(t)^{\dagger} \Xi(t) dU(t)
+ U(t)^{\dagger} d\Xi(t) dU(t) \nonumber\\
&=& U(t)^{\dagger} d\Xi(t) U(t)
+\Bigg\{ \Lambda(t) U(t)^{\dagger} dU(t)
+(\Lambda(t) U(t)^{\dagger} dU(t))^{\dagger}
\Bigg\} \nonumber\\
&& + \Bigg\{ U(t)^{\dagger} d\Xi(t) dU(t)+
(U(t)^{\dagger} d\Xi(t) dU(t))^{\dagger} \Bigg\} 
+ dU(t)^{\dagger} \Xi(t) dU(t).
\nonumber
\end{eqnarray}
Each term in the RHS is rewritten as follows:
\begin{eqnarray}
\Lambda(t)U(t)^{\dagger} dU(t)
&=&  \Lambda(t) \left( U(t)^{\dagger} dU(t) 
+\frac{1}{2} dU(t)^{\dagger} dU(t) \right) 
-\frac{1}{2} \Lambda(t) dU(t)^{\dagger} dU(t) \nonumber\\
&=& \Lambda(t) dA(t) 
+ \frac{1}{2} \Lambda(t) dA(t) dA(t),
\nonumber
\end{eqnarray}
where (\ref{eqn:antidA}) was used, and
\begin{eqnarray}
 U(t)^{\dagger} d\Xi(t) dU(t) 
&=& U(t)^{\dagger} d\Xi(t) U(t) dA(t), \nonumber\\
dU(t)^{\dagger} \Xi(t) dU(t) &=&
dA(t)^{\dagger} U(t)^{\dagger} \Xi(t) U(t) dA(t) 
= dA(t)^{\dagger} \Lambda(t) dA(t),
\nonumber
\end{eqnarray}
where (\ref{eqn:dU}) was used.
Then we have the equality 
\begin{eqnarray}
d \Lambda(t) &=& U(t)^{\dagger} d\Xi(t) U(t)
+ \Lambda(t) dA(t) + (\Lambda(t) dA(t))^{\dagger} \nonumber\\
&& + \frac{1}{2} \Lambda(t) dA(t) dA(t) 
+ \frac{1}{2} (\Lambda(t) dA(t) dA(t))^{\dagger} \nonumber\\
&& +U(t)^{\dagger} d\Xi(t) U(t) dA(t)
+(U(t)^{\dagger} d\Xi(t) U(t) dA(t))^{\dagger} 
+ dA(t)^{\dagger} \Lambda(t) dA(t).
\label{eqn:equality1}
\end{eqnarray}
The diagonal elements of (\ref{eqn:equality1}) give
\begin{eqnarray}
d \lambda_{i}(t) &=& \sum_{k, \ell} u_{ki}(t)^{*}
u_{\ell i}(t) d\xi_{k \ell}(t) \nonumber\\
&& + 2 \lambda_{i}(t) d \gamma_{ii}(t)+
d\phi_{ii}(t)+d\phi_{ii}(t)^{*}+d \psi_{ii}(t),
\quad 1 \leq i \leq N,
\label{eqn:equality2}
\end{eqnarray}
and the off-diagonal elements of (\ref{eqn:equality1}) give
\begin{eqnarray}
0 &=& \sum_{k, \ell} u_{k i}(t)^{*}u_{\ell j}(t) 
d\xi_{k \ell}(t)
+\lambda_{i}(t) d\alpha_{ij}(t) + \lambda_{j}(t) d \alpha_{ji}(t)^{*}
\nonumber\\
&&+ \lambda_{i}(t) d \gamma_{ij}(t)+\lambda_{j}(t) d \gamma_{ji}^{*}(t)
+ d \phi_{ij}(t)+ d\phi_{ji}(t)^{*} + d\psi_{ij}(t),
\quad 1 \leq i < j \leq N, \qquad
\label{eqn:equality3}
\end{eqnarray}
where we have used the notations
\begin{eqnarray}
d \gamma_{ij}(t) &\equiv& \left(\frac{1}{2} dA(t) dA(t) \right)_{ij}
= \frac{1}{2} \sum_{k} d \alpha_{ik}(t) 
d \alpha_{kj}(t) = d\gamma_{ji}(t)^{*},
\nonumber\\
\label{eqn:dphi1}
d \phi_{ij}(t) &\equiv& \Big( U(t)^{\dagger} d\Xi(t) U(t) dA(t)
\Big)_{ij} 
= \sum_{k, \ell, m} u_{k i}(t)^{*}
d\xi_{k \ell}(t) u_{\ell m}(t) d \alpha_{m j}(t), \\
d \psi_{ij}(t) &\equiv& \Big( dA(t)^{\dagger} \Lambda(t)
dA(t) \Big)_{ij} \nonumber\\
&=& \sum_{k} d \alpha_{ki}(t)^{*} \lambda_{k}(t) d \alpha_{kj}(t)
= - \sum_{k} d \alpha_{ik}(t) \lambda_{k}(t) d \alpha_{kj}(t).
\nonumber
\end{eqnarray}
Since $(\gamma_{ij}(t))_{1 \leq i, j \leq N},
(\phi_{ij}(t))_{1 \leq i, j \leq N},
(\psi_{ij}(t))_{1 \leq i, j \leq N}$ are functions
of finite variations, (\ref{eqn:equality1}) gives
\begin{eqnarray}
d \langle M_{i} \rangle_{t} &=&
\Bigg( U(t)^{\dagger} d \Xi(t) U(t) \Bigg)^{\dagger}_{ii}
\Bigg( U(t)^{\dagger} d \Xi(t) U(t) \Bigg)_{ii}
\nonumber\\
&=&
\sum_{k, \ell} \sum_{m, n}
u_{ki}(t)^{*} u_{\ell i}(t) u_{m i}(t)^{*} u_{n i}(t)
d\xi_{k \ell}(t) d \xi_{m n}(t).
\nonumber
\end{eqnarray}
On the other hand, for $dA(t)$ is anti-Hermitian, (\ref{eqn:equality3})
gives
\begin{eqnarray}
(\lambda_{j}(t)-\lambda_{i}(t)) d \alpha_{ij}(t)
&=& \sum_{k, \ell} u_{k i}(t)^{*} u_{\ell j}(t)
d\xi_{k \ell}(t) \nonumber\\
&+& (\lambda_{i}(t)+\lambda_{j}(t)) d\gamma_{ij}(t)
+d\phi_{ij}(t)+d\phi_{ji}(t)^{*}+d\psi_{ij}(t),
\label{eqn:equality4}
\end{eqnarray}
and using this equality we can rewrite (\ref{eqn:dphi1}) as
$$
d\phi_{ij}(t)= \sum_{k}(\lambda_{k}(t)-\lambda_{i}(t))
d \alpha_{ik}(t) d \alpha_{kj}(t).
$$
Then the finite-variation part of (\ref{eqn:equality2}) is 
written as
\begin{eqnarray}
&& 2 \lambda_{i}(t) d\gamma_{ii}(t)
+ d \phi_{ii}(t)+ d \phi_{ii}(t)^{*} + d\psi_{ii}(t) \nonumber\\
&=& \sum_{j} \Bigg\{ \lambda_{i}(t)
+2(\lambda_{j}(t)-\lambda_{i}(t))-\lambda_{j}(t) \Bigg\}
d \alpha_{ij}(t) d \alpha_{ji}(t) \nonumber\\
&=& \sum_{j} (\lambda_{j}(t)-\lambda_{i}(t))
d \alpha_{ij}(t) d \alpha_{ji}(t) \nonumber\\
&=& \sum_{j} \frac{1}{\lambda_{i}(t)-\lambda_{j}(t)}
{\bf 1}_{\{\lambda_{i}(t) \not= \lambda_{j}(t)\}}
\sum_{k, \ell, m, n}
u_{k i}(t)^{*} u_{\ell j}(t) u_{m j}(t)^{*}
u_{n i}(t) d \xi_{k \ell}(t) d \xi_{m n}(t),
\nonumber
\end{eqnarray}
where (\ref{eqn:equality4}) was used in the last equation.
This completes the proof. 
\hfill \qed

\vskip 0.5cm

\end{document}